\documentclass[12pt]{amsart} 
\usepackage{amsfonts, amssymb, amsmath, amsthm, 
latexsym, euscript, epic, eepic} 
\usepackage{amsrefs}
 \usepackage{times}

\author{Mark Holland}
\address{Mathematics Department, University of Manchester, Manchester 
UK}
\email{holland@maths.man.ac.uk}
\author{Stefano
Luzzatto}
\address{Mathematics Department, Imperial College, 
London SW7, UK}
\email{Stefano.Luzzatto@ic.ac.uk}
\urladdr{www.ma.ic.ac.uk/\textasciitilde luzzatto}

\begin{document}

\def \a{\alpha}
\newcommand{\g}{\gamma}
\def \l{\lambda}
\def \L{\Lambda}
\def \G{\Gamma}
\def \D{\Delta}
\def \d{\delta}
\def \e{\epsilon}
\def \g{\gamma}
\def \t{\theta}
\def \w{\omega}
\def \Om{\Omega}
\newcommand{\p}{\Phi}
\def \R{\mathbb{R}}
\def \Z{\mathbb{Z}}
\def \Boxx{\hfill$\blacksquare$}
\def \C{\EuScript{C}}
\def \E{\EuScript{E}}
\def \F{\EuScript{F}}
\def \I{\EuScript{I}}
\def \O{\EuScript{O}}
\def \A{\mathcal{A}}
\def \B{\mathcal{B}}
\def \K{\mathcal{K}}
\def \M{\mathcal{M}}
\def \N{\mathcal{N}}
\def \Q{\mathcal{Q}}
\def \T{\mathcal{T}}
\def \PP{\mathcal{P}}
\def \dd{\mathfrak{D}}

\theoremstyle{plain}
\newtheorem*{main theorem}{Main Theorem}
\newtheorem*{theorem*}{Theorem}
\newtheorem*{emptytheorem*}{}
\newtheorem{theorem}{Theorem}
\newtheorem{lemma}{Lemma}
\newtheorem{sublemma}{Sublemma}[lemma]
\newtheorem{proposition}{Proposition}[section]
\newtheorem*{proposition*}{Proposition}
\newtheorem{corollary}{Corollary}[section]
\newtheorem{claim}{Claim}[sublemma]

\theoremstyle{definition}
\newtheorem{definition}{Definition}
\newtheorem*{definition*}{Definition}
\newtheorem*{main definition}{Main Definition}
\theoremstyle{remark}
\newtheorem{remark}{Remark}
\newtheorem*{remark*}{Remark}
\newcommand{\mcal}[2]{\ensuremath{\mathcal #1^{(#2)}}}
\newcommand{\mcalp}[2]{\ensuremath{\mathcal #1^{(#2)}}(p)}

\newcommand{\pfi}[2]{\ensuremath{\partial_{#1}\Phi_{#2}}}
\newcommand{\pd}[2]{\ensuremath{\partial_{#1}^{#2}}}
\newcommand{\pdtil}[2]{\ensuremath{\tilde\partial_{#1}^{#2}}}

\title[A new proof of the Stable Manifold theorem]
{A new proof of the Stable Manifold Theorem \\ 
for hyperbolic fixed 
points on surfaces}
\date{\today}

\begin{abstract}
We introduce a new technique for proving the classical Stable Manifold 
theorem for hyperbolic fixed points. This method is much more 
geometrical than the standard approaches which rely on abstract fixed 
point theorems. It is based on the convergence of a canonical sequence 
of ``finite time local stable manifolds'' which are related to the 
dynamics of a finite number of iterations. 
\end{abstract}

\maketitle

\section{Introduction}
\subsection{Stable sets}
Let \( M \) be a two-dimensional
Riemannian manifold with Riemannian metric \( d \), and let \(
\varphi: M \to M \) be a   \( C^{2} \) diffeomorphism. 
Suppose that \( p\in M \) is a fixed point. 
A classical question concerns the effect that the existence of such a 
fixed point has on the global dynamics of \( f \). 
In particular we shall concentrate here on some properties of the set 
of  points which 
are forward asymptotic to \( p \).

\begin{definition}
The \emph{global stable set} \( W^{s}(p) \) of  \( 
p \) is
    \[ 
    W^{s}(z) = \{x\in M : 
    d(\varphi^{k} (x), 
    p) \to 0 \text{ as } k\to \infty\}.
    \]
\end{definition}

In general \( W^{s}(p) \) can be 
extremely complicated, both in its intrinsic geometry \cite{Kan94} 
and/or in the way it is embedded in \( M \) \cites{Poi90, Mos73, PalTak93}
It is useful therefore to begin with a study of that part of \( 
W^{s}(p) \)  which 
remain in a fixed neighbourhood of \( p \) for all forward 
iterations.  For \( \eta > 0 \) and \( k\geq 1 \) let 
   \[ 
   \mathcal N^{(k)}_{\eta}=\{x: d(\varphi^{j}(x), p) \leq 
   \eta \ \forall \ 0\leq j \leq k-1\}
\]
and
\[
   \mathcal N^{(\infty)}_{\eta}= \bigcap_{k\geq 1} \mathcal 
    N^{(k)}
   \]
\begin{definition}
        For  \( \eta > 0 \), 
define the \emph{local stable set} of  \( 
p \) by
 \[
    W^{s}_{\eta}(p) = W^{s}(p) \cap\mathcal 
    N^{(\infty)}_{\eta}(p)
 \]
    \end{definition}

In several situation it is possible to obtain a fairly comprehensive 
description of the geometrical and dynamical properties of the local 
stable manifold. In this paper we shall focus on the simplest setting 
of a hyperbolic fixed point. 
    We recall that the fixed point \( p \) is hyperbolic if the 
    derivative \( D\varphi_{p} \) has no eigenvalues on the unit circle. 
    
 \begin{theorem*}
     Let \( \varphi: M \to M \) be a \( C^{2} \) diffeomorphism of a 
     Riemannian surface and suppose that \( p \) is a hyperbolic 
     fixed point with eigenvalues 
\( 0<| \lambda_{s}|< 1 < | \lambda_{u}| \). 
 Then
     there exists a constant 
     \( \eta > 0 \) such that the following properties hold: 
     \begin{enumerate}
\item \( W^{s}_{\eta} (p) \) is \( C^{1+Lip} \) 
	 one-dimensional submanifold of M 
	tangent to \( E^{s}_{p} \);
\item \( |W^{s}_{\eta} (p)|\geq \eta \) 
	on either side of \( p \);
\item \( W^{s}_{\eta} (p) \) contracts at an 
	exponential rate.
	\item 
    \[ 
W^{s}_{\eta}(p) = \bigcap_{k\geq 0} \mathcal N^{(k)}_{\eta}(p). 
\]	
\end{enumerate}
     \end{theorem*}

  Perhaps the key statement here is that the local stable set is 
  actually a smooth submanifold of \( M \). 
  Notice that the global stable set can be written as a union of 
  preimages of \( W^{s}_{\eta} \) since any point whose orbit 
  converges to \( p \) must eventually remain in a neighbourhood of \( 
  p \) and therefore must eventually belong to \( W^{s}_{\eta} \). 
Thus we write
\[ 
W^{s}(p) = \bigcup_{n\geq 0} \varphi^{-n}(W^{s}_{\eta}(p)).
\]
The smoothness of \( W_{\eta}^{s} \) then implies that \( W^{s}(p) 
\) is also a smooth submanifold of \( M \). We remark that while the 
local stable manifold is an \emph{embedded} submanifold, i.e. a 
manifold in its own right in the induced topology, the global stable 
manifold is in general only an \emph{immersed} submanifold, i.e. a 
manifold in its intrinsic topology but not in the topology induced by 
the topology on \( M \). This is because  it may accumulate on itself 
and thus fail to be locally connected in the induced topology. 

\subsection{Main ideas}
The proof is based on the key notions of \emph{hyperbolic coordinates} 
and \emph{finite time local stable manifolds}. These are not standard 
concepts and therefore we describe them briefly here, leaving the 
details to the main body of the paper. The standard definition of 
hyperbolicity involves a decomposition of the tangent space into 
invariant subspaces, which coincide with the real \emph{eigenspaces} 
in the case of a fixed point. However this decomposition is 
specifically related to the \emph{asymptotic} properties of the dynamics and is not 
always the most useful or the most appropriate for studying the local 
geometry related to the dynamics for a finite number of iterations. 
Instead, elementary linear algebra arguments show that under 
extremely mild hyperbolicity 
conditions on the linear map \( D\varphi^{k}(z) \), 
there are well-defined \emph{most contracted} 
directions  \( e^{(k)}(z) \) which 
depend in a \( C^{1} \) manner on the point \( z \) if \( \varphi^{k} \) 
is \( C^{2} \). In particular,  they define a foliation of 
\emph{integral curves} \( \mathcal E^{(k)} \) which 
are \emph{most contracted curves under } \( \varphi^{k} \). 
They are the natural notion of local stable manifold relative to 
a finite number of iterations. 

In general, even for a hyperbolic fixed point, the finite time local 
stable manifolds will not coincide with the ``real'' asymptotic local 
stable manifold and will depend on the iterate \( k \). It is natural 
to expect however that the two concepts are related in the sense that 
finite time local stable manifolds converge to the asymptotic local 
stable manifold. In this paper we show that 
this is indeed  what happens. The argument does not require 
previous knowledge of the existence  of the asymptotic local stable 
manifold, and thus we obtain from the construction an alternative 
proof of the classical stable manifold theorem. 

The inspiration for this approach comes from the pioneering work of 
Benedicks and Carleson on non-uniformly hyperbolic dynamics in the 
 context of the H\'enon map \cite{BenCar91}, the 
generalization of this work to H\'enon-like maps by Mora 
and Viana \cite{MorVia93}, and the 
further refinement of these ideas in \cites{WanYou01, LuzVia}. 
In these papers, a sophisticated induction argument is developed to 
show that a certain non-uniformly hyperbolic structure is persistent, 
in a measure-theoretic sense, in typical parametrized families of maps. The 
induction necessarily requires some knowledge of 
geometrical structures based only 
a \emph{finite} number of iterations. 
Thus the notion of contractive directions is very natural and very 
efficient for producing dynamical foliations which incorporate and 
reflect information about the first \( k \) iterates of the map. 

We remark however that in the papers mentioned above, the 
construction of these  manifolds is very much 
embedded in the overall argument and there is no  focus on 
the intrinsic interest of this particular construction. Moreover, the 
specific characteristics of the maps under consideration, such as the 
strong dissipativity and specific distortion bounds, are used 
heavily in the estimates and it is not immediately clear  
what  precise assumptions are used in the construction.  
Thus the main purpose of the present paper can be seen as a first step 
in an attempt to draw on the ideas first 
introduced in \cite{BenCar91} and to formalize them and generalize them 
into a fully fledged theory of invariant  manifolds.

\subsection{Structure of the argument}

We start the proof with some relatively standard hyperbolicity and 
distortion estimates in a neighbourhood of \( p \) in Section 
\ref{local hyperbolicity}. 
In sections \ref{hyperbolic coordinates}  and \ref{finite time} 
we address the first main issue which is the convergence of finite 
time local stable manifolds. 
This will be addressed by carrying out several 
estimates which show that angle between successive most contracted 
direction \( e^{(k)} \) and \( e^{(k+1)} \) depends on the 
hyperbolicity of \( D\varphi^{k} \). If some hyperbolicity conditions 
are satisfied for all \( k \),  this sequence of angles 
forms a Cauchy sequence and therefore a limit direction \( e^{(\infty)} 
\) exists. We then show that we can also 
control the way in which contractive directions depend on the base 
point and eventually conclude that the sequence of most contracted 
leaves \( \mathcal E(k)(p) \) through the fixed point \( p \) 
actually converge (to something\ldots). 
In Section \ref{asymptotic} we begin to show that this limit ``object'' 
is indeed the local stable manifold of \( p \) by showing first of all 
that it has positive length. This 
requires a careful control of the size and geometry 
of the domains in which the 
contracting directions \( e^{(k)}(z) \) are defined, to ensure that 
 the finite time local stable manifold \( 
\mathcal E^{(k)} \) always have some fixed length. 
In section \ref{smoothness} we then show that this curve is 
sufficiently smooth. In section  \ref{s:contraction} we show that 
points converge exponentially fast to \( p \) and, finally, in
Section \ref{uniqueness} we prove uniqueness in the sense that the 
curve we have constructed is the only set of points which satisfy the 
properties of a local stable manifold. 
To simplify the notation we shall suppose that both 
eigenvalues are positive, 
the other cases are all dealt with in exactly the same way.

\subsection{Comparison with other approaches.}

The stable manifold theorem is one of the most basic results in the 
geometric theory of dynamical systems and differential equations. 
There exists several approaches and many generalizations 
\cites{Had01, Sma63, Irw70, Fen71, Fen74, Fen77, HirPugShu77, Pes76, 
Pes77, RueShu80,  
FatHerYoc83, Wig94, Rob95, Aba01, Cha01, AbbMaj01, Cha02}. 
We focus here on what we believe are the two main differences between 
our approach and existing ones. 

As far as 
we know, all existing approaches rely in one way or another on an 
application of the contraction mapping theorem to some suitable 
abstract space of \emph{candidate} local stable manifolds, formulated 
e.g. as the graphs or sequences with certain properties. A 
suitable complete metric and 
an operator depending  are then defined
in such a way that a fixed point under the operator corresponds to an 
invariant submanifold with the required properties. The existence and 
uniqueness of the fixed point are then a consequence of showing that 
the operator is a contraction and the application of the contraction 
mapping theorem. 

Our approach differs significantly from either of the approaches 
mentioned in at least two ways.  First of all we do not use the 
contraction mapping Theorem, not even in disguise. We show that the 
sequence of finite time local stable manifolds is Cauchy in an 
appropriate topology, by relating directly the distance between 
succesive leaves to the hyperbolicity. This is potentially a much 
more flexible approach and one which may be adaptable to situations 
with less uniform and perhaps subexponential forms of hyperbolicity, 
see comments in the next section. Secondly, the finite time local 
stable manifolds which approximate the real local stable manifold are 
canonically defined and have an important dynamical meaning in their 
own right. Again, this aspect of the construction might be useful 
in other situation and in applications.

\subsection{Generalizations and applications}

In this paper we present the construction of the local stable manifold 
in the simplest setting of a hyperbolic fixed point in dimension 2, in order 
to describe the ideas and the basic strategy in the clearest possible 
way. However we have no doubt that with some work, the basic argument 
will generalize in various directions. In work in progress, the 
authors are extending the techniques to cover the higher dimensional 
situation and the more general uniformly hyperbolic and projectively 
hyperbolic cases.  The smoothness results are not yet optimal since 
it is known that the stable manifold is as smooth as the map, however 
we believe that the optimal estimates can be recovered with more 
sophisticated higher order distortion estimates. It is also 
possible to study questions related to 
the dependence of the local stable manifolds on the base point or on 
some parameter by simply including the corresponding additional 
derivative estimates in the calculations, as already carried out for 
example in the papers mentioned above in which the original ideas of 
the integral curves of the most contracting directions was first 
introduced.

Perhaps the most important 
potential of this method, however, is to situations with very weak 
hyperbolicity.  Indeed, the convergence of the finite time local 
stable leaves is given by the Cauchy property which follows 
essentially 
from a summability condition on the derivative. 
Thus, in principle, this may be applicable along orbits for which the 
derivative does not necessarily admit exponential estimates. 
This fact is not completely explicit in the present proof since we 
are dealing with very hyperbolic example, but a future paper but the 
same authors will provide the minimal abstract conditions under which the 
convergence argument works. 

Finally we mention that there is also a vast amount of research 
literature concerning the explicit, often numerical, approximation 
of invariant manifolds.  We speculate that the approach given here 
might be more suitable than the classical theorems 
as a theoretical backdrop to these numerical studies, , 
The finite time local stable manifolds for example can probably be 
calculated with significant accuracy since they depend only a finite 
time information on the dynamics and the derivative. Our estimates 
then give conditions on how close these finite time manifolds are form 
the real thing.

\section{Local hyperbolicity and distortion estimates}
\label{local hyperbolicity}
The first step in the proof is to use the smoothness of \( \varphi \) 
to show that some hyperbolic structure exists in a 
neighbourhood of the point \( p \). 
First of all we introduce some notation which we will use extensively 
throughout the paper. 
 For \( z\in M \) and \( k\geq 1 \)
we let
\[ 
    F_{k}(z) = \|D\varphi^{k}_{z}\|\quad\text{ and } \quad E_{k}(z) =
    \|(D\varphi^{k}_{z})^{-1}\|^{-1} 
\]
denote the maximum expansion and the maximum contraction respectively 
of \(  D\varphi^{k}_{z}\). Then 
\[ H_{k}(z) =
    E_{k}(z) / F_{k}(z)
    \]
    denote the \emph{hyperbolicity} of \(  D\varphi^{k}_{z}\). 
  
    If a sequence of neighbourhoods \( \mathcal N^{(k)}  \) is fixed, 
    as it will be below, we shall often use the notation \( \bar 
    F_{k}, \bar E_{k} \) and \( \bar H_{k} \) to denote the maximum 
    values of these quantities in \( \mathcal H_{k} \). 
      Also, when expressing relationships between these quantities which hold 
    for all \( x \) in \( \mathcal N^{(k)} \) we shall often omit 
    explicit reference to \( x \).
    
    We shall also 
    use the same convention for other functions of \( x \) to be 
    introduced below. 
    We shall also use the notation 
    \[ 
    F_{j, k} = \|D\varphi^{k-j-1}_{\varphi^{j+1}(z)}\|,  
    \]
   for \( j=0, \ldots, k-1\)
The following Lemma 
follows from standard estimates in the theory of uniform 
hyperbolicity. We refer to \cite{KatHas94} for details and proofs. 

\begin{lemma}\label{regular growth}
    \(\exists \  K>0 \) such that 
\(\forall\   \delta  > 0 \),
\(\exists \  \varepsilon (\delta) > 0 \) 
such that for 
all \( k\geq j\geq 0 \) and all \( x\in \mathcal N^{(k)}_{\varepsilon} 
\) we have: 
\[ 
K ( \lambda_{u}+ \delta)^{j} \geq F_{j} \geq 
( \lambda_{u}- \delta)^{j} \geq 
(\lambda_{s}+ \delta)^{j} 
\geq  E_{j} \geq 
K^{-1} (\lambda_{s}- \delta)^{j}
\]
Moreover, we also have that 
\[ 
\sum_{j=0}^{k-1} F_{j} \leq K F_{k}; \quad
 F_{j} F_{j, k} \leq K F_{k}; \quad\text{ and } \quad
 \sum_{i =j}^{\infty} H_{i} \leq K H_{j} 
 \]
\end{lemma}
To simplify the notation below we shall generally use 
\[
 \lambda_{u}=\tilde\lambda_{u}-\tilde\varepsilon  
 \quad\text{and }\quad 
 \lambda_{s}=\tilde\lambda_{s}+\tilde\varepsilon. 
 \]
We will also some estimates on the higher order derivatives which 
follow easily from the hyperbolicity properties given above. 
\begin{lemma}\label{distortion}
 \( \exists \ K>0 \) such that   
 \(\forall \  x\in\mathcal N^{(k)} \) we have 
\[
\|D^{2}\varphi^{k}_{x}\| \leq K F^{2}_{k}; 
\quad \text{ and } \quad
\|D(\det
D\varphi^{k}_{x})\|
\leq K E_{k} F_{k}^{2}.  
\]

\end{lemma}

\begin{proof}
	Let \( A_{j}= D\varphi_{\varphi^{j}z} \) and let \( A^{(k)}=
	A_{k-1}A_{k-2}\dots A_{1} A_{0} \).  Let \( D A_{j}\) denote
	differentiation of \( A_{j} \) with respect to the space variables. 
	By the product rule for differentiation we have
\begin{equation}\label{productrule}
  \begin{split}
D^{2}\varphi^{k}_{z} 
 &
= DA^{(k)} 
= D(A_{k-1}A_{k-2}\dots
    A_{1}A_{0}
    ) 
  \\ &
    = \sum_{j=0}^{k-1} A_{k-1}\dots A_{j+1}(DA_{j}) A_{j-1}\dots A_{0}.
    \end{split}
\end{equation}
Taking norms on both sides of 
\eqref{productrule} and using the fact that 
\( A_{k-1}\dots A_{j+1} = D\varphi^{k-j-1}_{\varphi^{j+1}z} \), \(
A_{j-1}\dots A_{0} = D\varphi^{j}_{z} \) and, by the chain rule, \(
DA_{j} = D (D\varphi_{\varphi^{j}z}) = D^{2}\varphi_{\varphi^{j}z}
D\varphi^{j}_{z} \), we get
\begin{align*}
\|D^{2}\varphi^{k}_{z}\| 
&\leq \sum_{j=0}^{k-1}
\|D\varphi^{k-j-1}_{\varphi^{j+1}z}\| \cdot \|
D^{2}\varphi_{\varphi^{j}z}\| \cdot \|D\varphi^{j}_{z} \|^{2} 
\\
& \leq 
\sum_{j=0}^{k-1} \| D^{2}\varphi_{\varphi^{j}z}\| \hat F_{j} 
F_{j}^{2} 
\\
&\leq K \sum_{j=0}^{k-1}\hat F_{j} 
F_{j}^{2}.
\end{align*}
The last inequality follows from the fact that \(
\|D^{2}\varphi_{\varphi^{j}z}\| \) is uniformly bounded above.  
Then, by Lemma \ref{regular growth} we have 
\[ 
    \|D^{2}\varphi^{k}_{z}\| 
    \leq \sum_{j=0}^{k-1} \hat F_{j}F_{j}^{2} \leq 
    K\sum_{j=0}^{k-1} F_{k} F_{j} 
    \leq K F_{k} \sum_{j=0}^{k-1} F_{j} \leq K 
    F_{k}^{2}
    \]
    This proves the first part of the statement.  
To prove the second,  we argue along similar lines, this time letting \(
A_{j}= \det D\varphi_{\varphi^{j}z} \).  Then we have, as in
\eqref{productrule} above,
\(
D(\det\varphi^{k}_{z}) = DA^{(k)} = \sum_{j=0}^{k-1} A_{k-1} \dots
A_{j+1}(DA_{j}) A_{j-1}\dots A_{0}.  \)
Moreover we have that
\( A_{k-1}\dots A_{j+1} = \det D\varphi^{k-j-1}_{\varphi^{j+1}z} \),
\( A_{j-1}\dots A_{0} = \det D\varphi^{j}_{z} \) and, by the chain
rule, also \( DA_{j} = D(\det D\varphi_{\varphi^{j}z}) = (D \det
D\varphi_{\varphi^{j}z}) D\varphi^{j}_{z} \). This gives
\begin{equation}\label{det}
    D(\det\varphi^{k}_{z}) = \sum_{j=0}^{k-1} (\det
    D\varphi^{k-j-1}_{\varphi^{j+1}z}) (D \det D\varphi_{\varphi^{j}z})
    (\det D\varphi^{j}_{z} )(D\varphi^{j}_{z}).
\end{equation}
By the multiplicative property of the determinant we have
the equality \[
(\det D\varphi^{k-j-1}_{\varphi^{j+1}z}) (\det D\varphi^{j}_{z} ) =
\det D\varphi^{k}_{z}/\det D\varphi_{\varphi^{j}{z}}.  
\] 
Therefore we get 
\begin{align*} 
\| D(\det D\varphi^{k}_{z})\| 
&\leq \ |\det D\varphi^{k}_{z} |
\sum_{j=0}^{k-1} \frac{\|D(\det D\varphi_{\varphi^{j}(z)}\|}{|\det
D\varphi_{\varphi^{j}(z)}|} F_{j}
\\
&\leq  E_{k}F_{k} \sum_{j=0}^{k-1} F_{j} 
\leq K E_{k} F_{k}^{2}.
\end{align*}
The first inequality follows by taking
norms on both sides of \eqref{det}, 
the second inequality follows from the fact that \( \|D(\det
D\varphi_{\varphi^{j}(z)}\| \) and \( |\det D\varphi_{\varphi^{j}(z)}|
\) are uniformly bounded above and below and the fact that 
that \( \det D\varphi^{k} = E_{k} F_{k}
\), the third inequality follows from an application of Lemma 
\ref{regular growth}. 
\end{proof}

\section{Hyperbolic coordinates}\label{hyperbolic coordinates}
In the context of hyperbolic fixed points (or general uniformly 
hyperbolic sets) we are used to thinking of the eigenspaces (or the 
subspaces given by the hyperbolic decomposition) as providing the 
basic axes or coordinate system associated to the hyperbolicity. 
However this is not necessarily the most natural splitting of the space. 
Indeed, the hyperbolicity condition 
\(
H_{k} = E_{k}/F_{k} < 1
\)
(notice that we always have \( H_{k}\leq 1 \)) 
implies that the linear map \( D\varphi^{k} \) maps 
the unit circle \(\mathcal S \subset T_{z}M \) 
to an (non-circular) ellipse \( \mathcal S_{k} =
D\varphi_{z}^{k}(\mathcal S) \subset T_{\varphi^{k}(z)}M \)
with well defined major and minor axes. 
    \begin{figure}[ht]
\begin{center}
\setlength{\unitlength}{0.00038333in}
\begingroup\makeatletter\ifx\SetFigFont\undefined%
\gdef\SetFigFont#1#2#3#4#5{%
  \reset@font\fontsize{#1}{#2pt}%
  \fontfamily{#3}\fontseries{#4}\fontshape{#5}%
  \selectfont}%
\fi\endgroup%
{\renewcommand{\dashlinestretch}{30}
\begin{picture}(9313,2842)(0,-10)
\put(1500,1282){\ellipse{2424}{2424}}
\path(6000,82)(8100,2557)(6000,82)
\path(6675,1732)(7500,907)(6675,1732)
\path(1050,2407)(2025,157)
\path(450,682)(2550,1807)(2625,1807)
\drawline(2625,1882)(2625,1882)
\path(3225,1432)(3227,1433)(3233,1435)
	(3242,1439)(3256,1445)(3276,1453)
	(3299,1463)(3327,1474)(3358,1486)
	(3391,1499)(3425,1513)(3458,1526)
	(3491,1538)(3523,1550)(3553,1562)
	(3582,1572)(3609,1582)(3635,1591)
	(3660,1599)(3683,1606)(3707,1613)
	(3729,1620)(3752,1626)(3775,1632)
	(3798,1638)(3822,1643)(3846,1648)
	(3870,1653)(3896,1658)(3922,1663)
	(3948,1668)(3976,1672)(4004,1676)
	(4033,1680)(4062,1684)(4091,1688)
	(4120,1691)(4150,1694)(4179,1696)
	(4208,1699)(4237,1701)(4265,1702)
	(4293,1704)(4320,1705)(4347,1706)
	(4373,1707)(4399,1707)(4425,1707)
	(4451,1707)(4477,1707)(4503,1706)
	(4530,1705)(4557,1704)(4585,1702)
	(4613,1701)(4642,1699)(4672,1696)
	(4701,1694)(4731,1691)(4761,1688)
	(4791,1684)(4820,1680)(4849,1676)
	(4878,1672)(4906,1668)(4934,1663)
	(4961,1658)(4987,1653)(5013,1648)
	(5038,1643)(5063,1638)(5088,1632)
	(5114,1626)(5141,1619)(5167,1612)
	(5194,1604)(5222,1597)(5249,1589)
	(5277,1580)(5305,1572)(5333,1563)
	(5361,1554)(5388,1544)(5415,1535)
	(5440,1526)(5466,1517)(5490,1509)
	(5512,1500)(5534,1492)(5555,1485)
	(5574,1477)(5592,1470)(5609,1463)
	(5625,1457)(5649,1448)(5670,1439)
	(5691,1430)(5711,1421)(5731,1413)
	(5751,1403)(5772,1394)(5792,1385)
	(5811,1376)(5850,1357)
\path(5728.982,1382.587)(5850.000,1357.000)(5755.260,1436.526)
\path(6000,82)(6001,81)(6004,78)
	(6013,71)(6025,62)(6039,51)
	(6054,42)(6068,34)(6082,28)
	(6096,23)(6113,19)(6124,17)
	(6136,16)(6149,14)(6163,13)
	(6179,12)(6196,12)(6213,12)
	(6232,12)(6251,13)(6271,14)
	(6290,16)(6310,19)(6330,21)
	(6349,24)(6368,28)(6388,32)
	(6403,36)(6419,40)(6435,44)
	(6452,49)(6470,55)(6488,61)
	(6507,68)(6527,75)(6547,84)
	(6568,92)(6589,102)(6610,112)
	(6631,122)(6652,133)(6672,144)
	(6693,156)(6714,168)(6734,181)
	(6755,194)(6775,207)(6792,219)
	(6810,231)(6828,244)(6847,257)
	(6866,271)(6885,286)(6905,301)
	(6925,317)(6946,334)(6967,351)
	(6988,368)(7009,385)(7030,403)
	(7050,421)(7071,438)(7090,456)
	(7110,473)(7128,490)(7146,506)
	(7163,522)(7180,538)(7195,553)
	(7210,568)(7225,582)(7242,599)
	(7258,616)(7274,633)(7290,650)
	(7306,668)(7322,686)(7339,706)
	(7356,727)(7375,749)(7393,772)
	(7412,795)(7431,819)(7448,841)
	(7464,861)(7477,878)(7487,891)
	(7494,900)(7498,905)(7500,907)
\path(6000,82)(5999,85)(5998,92)
	(5995,103)(5992,120)(5988,141)
	(5983,165)(5978,191)(5974,217)
	(5970,243)(5966,267)(5964,290)
	(5962,311)(5960,331)(5960,351)
	(5960,369)(5961,388)(5963,407)
	(5964,424)(5967,442)(5970,461)
	(5974,480)(5978,500)(5983,521)
	(5988,543)(5995,566)(6001,589)
	(6009,612)(6016,636)(6025,660)
	(6033,684)(6042,708)(6051,731)
	(6061,754)(6070,777)(6080,800)
	(6090,822)(6100,844)(6110,865)
	(6119,886)(6130,907)(6140,928)
	(6152,951)(6163,973)(6176,997)
	(6188,1020)(6201,1044)(6215,1068)
	(6228,1093)(6242,1117)(6256,1141)
	(6270,1165)(6284,1188)(6297,1210)
	(6311,1232)(6324,1253)(6337,1274)
	(6350,1294)(6362,1313)(6375,1332)
	(6388,1351)(6400,1369)(6413,1388)
	(6426,1407)(6440,1425)(6454,1445)
	(6468,1464)(6483,1483)(6499,1503)
	(6514,1523)(6530,1543)(6546,1562)
	(6563,1582)(6579,1601)(6596,1620)
	(6612,1639)(6629,1657)(6646,1675)
	(6662,1693)(6679,1710)(6695,1727)
	(6713,1744)(6728,1760)(6745,1776)
	(6762,1793)(6779,1809)(6797,1826)
	(6816,1844)(6835,1862)(6855,1880)
	(6875,1898)(6896,1917)(6917,1935)
	(6937,1954)(6958,1972)(6979,1990)
	(7000,2008)(7020,2025)(7040,2042)
	(7059,2058)(7078,2074)(7096,2089)
	(7113,2104)(7130,2118)(7147,2131)
	(7163,2144)(7181,2160)(7200,2175)
	(7218,2189)(7236,2204)(7254,2218)
	(7273,2233)(7291,2247)(7309,2261)
	(7327,2274)(7345,2288)(7362,2301)
	(7380,2313)(7396,2325)(7412,2336)
	(7428,2347)(7443,2358)(7458,2367)
	(7472,2377)(7486,2386)(7500,2394)
	(7517,2405)(7535,2416)(7553,2426)
	(7571,2437)(7590,2447)(7608,2457)
	(7628,2467)(7647,2477)(7665,2486)
	(7683,2495)(7701,2503)(7717,2510)
	(7733,2516)(7748,2522)(7762,2527)
	(7775,2532)(7792,2537)(7808,2542)
	(7824,2546)(7840,2549)(7856,2552)
	(7872,2554)(7886,2555)(7901,2556)
	(7914,2557)(7927,2557)(7938,2557)
	(7950,2557)(7962,2557)(7973,2557)
	(7986,2557)(8001,2557)(8017,2557)
	(8035,2557)(8053,2557)(8071,2557)
	(8086,2557)(8095,2557)(8099,2557)(8100,2557)
\path(8100,2557)(8101,2556)(8104,2553)
	(8111,2544)(8121,2532)(8133,2518)
	(8143,2503)(8153,2489)(8161,2475)
	(8168,2461)(8175,2444)(8179,2433)
	(8184,2421)(8188,2408)(8192,2393)
	(8196,2377)(8200,2360)(8203,2341)
	(8206,2322)(8208,2301)(8209,2280)
	(8210,2258)(8210,2236)(8209,2213)
	(8207,2191)(8204,2168)(8200,2144)
	(8197,2127)(8193,2109)(8188,2090)
	(8183,2071)(8176,2050)(8170,2029)
	(8162,2007)(8153,1983)(8144,1959)
	(8134,1935)(8124,1909)(8112,1883)
	(8100,1857)(8088,1830)(8075,1804)
	(8062,1777)(8048,1751)(8034,1725)
	(8020,1698)(8005,1672)(7990,1646)
	(7975,1619)(7963,1598)(7950,1577)
	(7936,1555)(7922,1532)(7907,1509)
	(7891,1484)(7874,1458)(7856,1430)
	(7837,1401)(7816,1370)(7794,1337)
	(7771,1302)(7746,1265)(7720,1227)
	(7694,1189)(7667,1149)(7640,1110)
	(7614,1073)(7590,1038)(7568,1005)
	(7549,977)(7533,954)(7520,936)
	(7511,922)(7505,914)(7501,909)(7500,907)
\put(3900,1882){\makebox(0,0)[lb]{\smash{{{\SetFigFont{6}{7.2}{\rmdefault}{\mddefault}{\updefault}$D\p^{k}(z)$}}}}}
\put(825,2557){\makebox(0,0)[lb]{\smash{{{\SetFigFont{6}{7.2}{\rmdefault}{\mddefault}{\updefault}$e^{(k)}$}}}}}
\put(0,382){\makebox(0,0)[lb]{\smash{{{\SetFigFont{6}{7.2}{\rmdefault}{\mddefault}{\updefault}$f^{(k)}$}}}}}
\put(8175,2632){\makebox(0,0)[lb]{\smash{{{\SetFigFont{6}{7.2}{\rmdefault}{\mddefault}{\updefault}$f^{(k)}_{k}$}}}}}
\put(7650,682){\makebox(0,0)[lb]{\smash{{{\SetFigFont{6}{7.2}{\rmdefault}{\mddefault}{\updefault}$e^{(k)}_{k}$}}}}}
\put(1725,1132){\makebox(0,0)[lb]{\smash{{{\SetFigFont{6}{7.2}{\rmdefault}{\mddefault}{\updefault}$z$}}}}}
\put(7200,1282){\makebox(0,0)[lb]{\smash{{{\SetFigFont{6}{7.2}{\rmdefault}{\mddefault}{\updefault}$\p^kz$}}}}}
\end{picture}
}
\end{center}
\end{figure}
The unit vectors \( e^{(k)}, f^{(k)} \) 
which are mapped to the minor and major
axis respectively of the ellipse, and are thus the \emph{most 
contracted} and \emph{most expanded} vectors respectively,  are given
analytically as solutions to
    the differential equation \( d\|D\varphi_{z}^{k}(\cos\theta, 
    \sin\theta)\|/d\theta = 0 
    \) which can be solved to give the explicit formula
    \begin{equation}\label{contractive directions}
\tan 2\theta  =
\frac{2 (\pfi x1^{k}\pfi y1^{k} +\pfi x2^{k}\pfi y2^{k})}
{(\pfi x1^{k})^2+(\pfi x2^{k})^2 - (\pfi y1^{k})^2 -(\pfi y2^{k})^2}.
\end{equation}
In particular, \( e^{(k)} \) and \( f^{(k)} \) are always \emph{orthogonal}
 and clearly \emph{do not} in general correspond to the 
stable and unstable eigenspaces of \( D\varphi^{k} \). 
We shall adopt the notation 
\[      
e^{(k)}_{j}= D\varphi^{j}(e^{(k)}) 
\quad\text{ and }
f^{(k)}_{j} = D\varphi^{j}(e^{(k)})
\]
where everything is of course relative to some given point \( x \). 
Notice that 
\[ 
E_{k} = \|(D\varphi^{k}_{z})^{-1}\|^{-1} = 
\|e^{(k)}_{k}(z)\| = \|D\varphi^{k}_{z}(e^{(k)}(z)) \|
\]
and 
\[
F_{k}(z) = \|D\varphi^{k}_{z}\| = \|f^{(k)}_{k}(z)\| =
\|D\varphi^{k}_{z}(f^{(k)}(z)) \|.
\]
Thus \( H_{k} \)
is a natural way of expressing the overall \emph{hyperbolicity} of \(
\varphi^{k} \) at \( z \).

We shall use \( K \) to denote a generic constant which is
allowed to depend only on the diffeomorphism \( \varphi \). 
For simplicity we shall allow this to change even within one 
inequality when this doe not create any ambiguities or confusion. 
We also define angles \( \phi^{(k)},  \theta^{(k)} \) by
\[
\phi^{(k)}=\measuredangle(e^{(k)},e^{(k+1)}),  \quad
\quad
\phi^{(k)}_j=\measuredangle(e^{(k)}_j,e^{(k+1)}_j)
\]
and
\begin{gather*}
e^{(k)}=(\cos\t^{(k)},\,\sin\t^{(k)}),\quad 
f^{(k)}=(-\sin\t^{(k)},\,\cos\t^{(k)}).\\
e^{(k)}_{k}=E_k (\cos\t^{(k)}_{k},\,\sin\t^{(k)}_{k}),\quad 
f^{(k)}_{k}=F_k (-\sin\t^{(k)}_{k},\,\cos\t^{(k)}_{k}).\\
\end{gather*}
Notice that 
\[ 
\theta^{(k)}= \theta^{(1)}+\sum_{j=1}^{k-1}\phi^{(j)}.
\]
and 
\[ 
De^{(k)}= (-D\theta^{(k)} (\sin \theta^{(k)}), 
D\theta^{(k)}(\cos\theta^{(k)}))
\]
which implies in particular that 
\begin{equation}\label{norms}
    \|De^{(k)}\|\leq K \|D\theta^{(k)}\| \leq K \|D\theta^{(1)}\|+ 
    K \sum_{j=1}^{k-1} \|D\phi^{(j)}\|.
\end{equation} 
where the constant \( K \) depends only on the choice of the norms.


The notion of most contracted and most expanded directions are 
central to the  approach to the stable manifold theorem given 
here.  
 We start with a lemma concerning the convergence 
 of the sequence of most 
 contracted direction \( e^{(k)}(x) \) as \( k \to \infty \).

\begin{lemma} \label{angleconvergence}
\(\exists \  K > 0 \) 
such that  \( \forall \ k\geq 1 \) and
\( x\in \mathcal N^{(k+1)} \) and \( \phi^{(k)}=\phi^{(k)}(x)
 \), we have
 \[
 |\phi^{(k)}| \leq K H_{k}
\]
\end{lemma}

\begin{proof}
Write $e^{(k)}=\eta e^{(k+1)}+\varphi f^{(k+1)}$ where 
$\eta^2+\varphi^2=1$ by normalization.
Linearity implies that $e^{(k)}_{k+1}=
\eta e^{(k+1)}_{k+1}+\varphi f^{(k+1)}_{k+1}$ and orthogonality implies that
\(
\|e^{(k)}_{k+1}\|^2 =\eta^2\|e^{(k+1)}_{k+1}\|^2+\varphi^2
\|f^{(k+1)}_{k+1}\|^2 
 =\eta^2 E^{2}_{k+1}+\varphi^2 F^{2}_{k+1} 
\)
where $E_k=\|e^{(k)}_k\|,\,F_k=\|f^{(k)}_k\|$.
Since $\phi^{(k)}=\tan^{-1}(\varphi/\eta)$ we get 
$$
|\tan\phi^{(k)}|=
\left(\frac{\|e^{(k)}_{k+1}\|^2-E^{2}_{k+1}}
{F^{2}_{k+1}-\|e^{(k)}_{k+1}\|^2}\right)^{\frac{1}{2}} \leq 
\frac{\|e^{(k)}_{k+1}\|/F_{k+1}}
{\left(1-\|e^{(k)}_{k+1}\|^2/F^{2}_{k+1}\right)^{\frac{1}{2}}}
$$
now notice that we can choose \( k_{0} 
 \)  large enough depending only on the map, 
 so that 
 \(  \|D\varphi_{\varphi^{k}(x)}\| E_{k}/F_{k+1}< 1/2\)
 for all \( k\geq k_{0} \). Then, for \( k \geq k_{0} \) we have 
 \[ 
 |\tan\phi^{(k)}|= \leq K \|e^{(k)}_{k+1}\|/F_{k+1} \leq KH_{k}
 \]
 where we have used the boundedness of the derivative to obtain the 
 last equality. Since \( k_{0}  \) depends only on \( \varphi \), we 
 can take care of all iterates \( k\leq k_{0} \) by simply adjusting 
 the constant \( K \). 
\end{proof}

The next result says gives some control over the dependence of the 
most contracting directions on the base point. 

\begin{lemma}\label{derivativeconvergence}
 \( \exists \ K> 0 \) such that  \( \forall \ k\geq 1 \) 
\(\forall\  x\in \mathcal N^{(k+1)} \) and \(
    \phi^{(k)}=\phi^{(k)}(x) \), we have
\[
\|D\phi^{(k)}\| \leq K E_{k}\quad\text{and}\quad \|De^{(k)}\| \leq K 
\]
\end{lemma}
 
The proof of Lemma \ref{derivativeconvergence} is relatively technical 
and we postpone it to the final Section. First though we prove another 
estimate which will be used below.

\begin{lemma}\label{anglecontraction}
   \( \exists \ K>0 \) such that 
     \(\forall\  x\in \mathcal N^{(k)} \) and \(\forall\   0\leq j \leq 
    k \) we have 
\[ 
\|e^{(k)}_{j}\|\leq K F_{j} \sum_{i=j}^{k-1} H_{i} + E_{k} \leq K E_{j}.
\]
    \end{lemma}

   \newcommand{\coni}[1]{\ensuremath{e^{(#1)} }}
 \newcommand{\exi}[1]{\ensuremath{f^{(#1)} }}
 \newcommand{\con}[2]{\ensuremath{e^{(#1)}(#2) }}
 \newcommand{\ex}[2]{\ensuremath{f^{(#1)}(#2) }}
 \newcommand{\conim}[3]{\ensuremath{e^{(#1)}_{#2}(#3) }}
 \newcommand{\exim}[3]{\ensuremath{f^{(#1)}_{#2}(#3) }}
 \newcommand{\conimi}[2]{\ensuremath{e^{(#1)}_{#2} }}
 \newcommand{\eximi}[2]{\ensuremath{f^{(#1)}_{#2} }}

\begin{proof}
   Using the linearity of the derivative we write
    \[
    \|e^{(k)}_{j}\|=
    \|D\varphi^{j} \coni k\| 
    \leq
\|D\varphi^{j} (\coni k-\coni j)\| + \|D\varphi^{j} \coni j\|. 
\] 
Then we have 
\( \| D\varphi^j \coni j\|= \|\conimi jj\| = E_{j}\) and 
\[
\|D\varphi^j (\coni
k-\coni j)\|\leq \|D\varphi^j \| \|\coni k-\coni j\| \leq
F_{j} \|\coni k-\coni j\|.
\]  
By Lemma \ref{angleconvergence} we have 
\[\|\coni k-\coni j\|\leq \sum_{i=j}^{k-1} \phi^{(i)} 
\leq K \sum_{i=j}^{k-1} H_{i}.
\]
This give the first inequality.  The second follows from 
Lemma \ref{regular growth}.
\end{proof}

\section{Finite time local stable manifolds}
\label{finite time}
The smoothness of the field of most contracted directions implies, by 
standard results on the existence of local integrals for vector 
fields, the existence of \emph{most contracted} and \emph{most 
expanded} leaves. Thus for every \( k\geq 0 \) we consider the 
integral curve to the field \( e^{(k)} \) of contractive directions 
which are well defined in \( \mathcal N^{(k)} \). 
\begin{definition}
\( \mathcal E^{(k)}(p) \)
    is called the \emph{time \( k \) local stable manifold of  of \( p\)}.
    \end{definition}
These \emph{finite time} local stable manifolds have not been used 
much in dynamics even though from a certain point of view they are 
even more natural than the traditional local stable manifolds. Indeed, 
the local stable manifolds captures information about some local 
geometrical structure related to the asymptotic nature of the 
dynamics, while perhaps not containing the right kind of information 
regarding geometrical structures associated to some finite number of 
iterations of the map. As far as the dynamics for some finite number 
of iterates is concerned, the finite time local stable manifold \( 
\mathcal E^{(k)} \) is the natural object to look at, since it is a 
canonically defined submanifold which is locally \emph{most 
contracted} under a finite number of iterations of the map.

We shall use the notation \(  \mathcal E^{(k)}_{\pm}(p) \) to denote 
the two pieces of \( \mathcal E^{(k)}(p) \) on either side of \( p\). 
\emph{A priori} we have no means of saying how successive leaves \( 
\mathcal E^{(k)} \) are related nor, importantly, what their length 
are. Indeed, the neighbourhoods \( \mathcal N^{(k)} \) are in 
principle shrinking in size and thus the length of the leaves \( 
\mathcal E^{(k)} \) could do the same. 
The heart of our argument however is precisely to show that this does 
not happen. First we  apply the pointwise convergence results 
of the previous section to show that the sequence of curves \( 
\{\mathcal E^{(k)}(p)\} \) is a uniformly convergent Cauchy sequence 
of curves and thus in particular converges pointwise to a limit curve 
\( \mathcal E^{(\infty)}(p) \). In Section \ref{asymptotic} 
we show that the arclength of 
the curves \( 
\{\mathcal E^{(k)}(p)\} \)  on both sides of \( p \) 
is uniformly bounded below, implying that 
the limit curve \( \mathcal E^{(\infty)}(p) \) has positive length. 
Then, in the final Sections we complete the proof by showing that 
\( \mathcal E^{(\infty)}(p) \) is the local stable manifold of \( p 
\) and has the required properties.

Let $z^{(k)}_t$ and $z^{(k+1)}_t$, be parametrizations by arclength of
the two curves \( \mathcal E^{(k)}(p) \) and \( \mathcal E^{(k+1)}(p)
\) respectively, with \( z^{(k)}_{0}= z^{(k+1)}_{0} = p \) and choose 
\( t_{0} \) so that both \( \{z^{(k+1)}_{t}\}_{t=-t_{0}}^{t_{0}}\) 
and \( \{z^{(k)}_{t}\}_{t=-t_{0}}^{t_{0}}\) are both contained in 
\( \mathcal N^{(k+1)} \). We shall prove below that we can choose \( 
t_{0} \) uniformly in \( k \). For the  moment we obtain an estimate 
for the distance between the two curves. 
\begin{lemma}
     \label{leafcontraction}
\(\exists \  K>0 \) such that 
 \(\forall \  k\geq 1 \) and \( -t_{0}\leq t \leq t_{0}
\) we have
    \[ 
    |z^{(k)}_{t}-z^{(k+1)}_{t}| \leq  
    Kt \bar H_{k} e^{t K}\leq Kt_{0} \bar H_{k}e^{t_{0} K}  
    \] 
\end{lemma}
\begin{proof}
Writing
      \(
      z^{(k)}_{t} = z + \int_{0}^{t}e^{(k)}(z_{s}) ds \) and \( z^{(k+1)}_t
      = z  + \int_{0}^{t}e^{(k+1)}( z^{(k_+1)}_{s}) ds \) 
      we have 
      \(   |z^{(k)}_{t} - z^{(k+1)}_{t}| = \int_{0}^{t}
\|e^{(k)}(z^{(k)}_{s}) -e^{(k+1)}(z^{(k+1)}_{s}) \| ds \). Thus, by 
the triangle inequality, this gives
      \begin{equation}\label{integral}
     \begin{split}  |z^{(k)}_{t} - z^{(k+1)}_{t}| 
     \leq   
      \int_{0}^{t} &
    \|e^{(k)}(z^{(k)}_{s}) -e^{(k)}(z^{(k+1)}_{s})\| \\ 
    + &
    \|e^{(k)}(z^{(k+1)}_{s})
     -e^{(k+1)}(z^{(k+1)}_{s})\| ds.
     \end{split}
 \end{equation}
  By the Mean Value Theorem and Lemma \ref{derivativeconvergence} we 
     have
\begin{equation}\label{eq1}
    \begin{split}
    \|e^{(k)}(z^{(k)}_{s}) -e^{(k)}(z^{(k+1)}_{s})\|  
    &\leq  K
    \|De^{(k)}\|  \ |z^{(k)}_{s}-z^{(k+1)}_{s}| 
    \\
    &\leq K 
    |z^{(k)}_{s}-z^{(k+1)}_{s}|
    \end{split}
\end{equation}
    and, by Lemma \ref{angleconvergence} we have
\begin{equation}\label{eq2}
\|e^{(k)}(z^{(k+1)}_{s})
     -e^{(k+1)}(z^{(k+1)}_{s})\| = |\phi^{(k)}| \leq K \bar H_{k}.
\end{equation}
From  \eqref{eq1} and \eqref{eq2} we get 
\begin{equation}\label{uniformconvergence}
    \|e^{(k)}(z_{s}^{(k)}) - e^{(k+1)}(z_{s}^{(k+1)})\|
    \leq K (|z_{s}^{(k)}-z_{s}^{(k+1)}| + H_{k}).
\end{equation}
Substituting \eqref{uniformconvergence} 
into \eqref{integral} 
and using Gronwall's inequality gives
 \begin{equation}\label{uniformpointwiseconvergence}
|z^{(k)}_{t}- z^{(k+1)}_{t}| \leq K \int_{0}^{t} 
    |z^{(k)}_{s}-z^{(k+1)}_{s}| + \bar H_{k} ds  
    \leq 
    K t \bar H_{k} e^{t K} 
\end{equation}
\end{proof}

\section{The asymptotic local stable manifold}
\label{asymptotic}
Lemma \ref{leafcontraction} shows that the finite time local stable 
manifolds are exponentially close as long as they are both of some 
positive length. It does not however imply that this length can be 
guaranteed. To show this we fix first of all a sequence \( \omega_{k} 
\) where 
\[ 
\omega_{k} = Ke^{\eta K} \eta \bar H_{k}
\]
where the constant \( K \) is chosen as in Lemma \ref{leafcontraction}. 
Then we let  \( \mathcal 
T_{\omega_{k}}(\mathcal E^{(k)}(p)) \) denote the neighbourhood of 
\( \mathcal E^{(k)}(p) \) of size \( \omega_{k} \).  We now set up and 
inductive argument. First of all we define the inductive condition
\begin{equation}\tag*{\( (*)_{k} \)}
    |\mathcal E^{(k)}_{\pm}(p)| \geq \eta \quad\text{ and }
\quad \mathcal T_{\omega_{k}}(\mathcal E^{(k)}(p)) \subset 
\mathcal N^{(k+1)}
\end{equation}
Notice that the contractive directions  \( e^{(1)} \) are well defined in 
\( \mathcal N^{(1)}=B_{\varepsilon}(p) \) by Lemma \ref{regular 
growth}. Therefore the time 1 local stable manifold \( \mathcal 
E^{(1)}(p) \) exists and, by taking \( \eta \) sufficiently small, we 
can guarantee that condition \( (*)_{1} \) is satisfied. Thus it just 
remain to prove the general inductive step

\begin{lemma}
\( \forall \ \eta>0 \) sufficiently small and \( k\geq 1 \), 
\[ 
(*)_{k} \ \Rightarrow \ (*)_{k+1}.
\]
\end{lemma}

\begin{proof}
    We assume condition \( (*)_{k} \) and start by proving 
    that \( \mathcal E^{(k+1)}(p) \) exists and has the required
length.  
By the assumption that 
\( \mathcal T_{\omega_{k}}(\mathcal E^{(k)}(p)) \subset 
\mathcal N^{(k+1)} \),
the vector field of contractive direction 
\( e^{(k+1)} \) is well defined in the whole of 
\( \mathcal T_{\omega_{k}}(\mathcal E^{(k)}(p))  \). 
Therefore \(  (\mathcal E^{(k+1)}(p)) \) certainly exists and the 
question is whether it satisfies the condition on the length. 
The fact that is does follows by Lemma 
\ref{leafcontraction} and our choice of the sequence \( \omega_{k} \). 
Indeed, Lemma \ref{leafcontraction} implies that 
    \( z^{(k+1)}_{t} \in \mathcal T_{\omega_{k}}(\mathcal E^{(k)}(p))  \)
    for all \( t\leq t_{0} \) where \( t_{0} \) can be chosen as \( 
    t_{0}= \eta \). 
    
    To show the second part of the statement, it is 
    sufficient to show that 
    \(x\in  \mathcal T_{\omega_{k+1}}(\mathcal 
    E^{(k+1)}(p)) \) implies \( x\in\mathcal N^{(k+2)} 
    \)
i.e.
\[
    d(p, \varphi^{j}(x)) 
    \leq \varepsilon \quad \forall\  j\leq k+1 
    \]
    First of all, since \( x \in  \mathcal T_{\omega_{k+1}}(\mathcal 
    E^{(k+1)}(p)) \), we can 
    choose some \( y \in \mathcal E^{(k+1)}(p) \) such that 
    \( d(x,y) \leq \omega_{k+1} \). Then, 
    by the triangle inequality 
    we have 
       \begin{equation}\label{triangle}
       d(p, \varphi^{j}(x)) 
   \leq d(p, \varphi^{j}(y)) + d(\varphi^{j}(y), 
   \varphi^{j}(x)).  \end{equation}
In the 
   following inequalities we always assume that all values depending 
   on \( x \) are taken to be the maximum over all \( x \) in \( 
   \mathcal T_{\omega_{k+1}} (\mathcal 
    E^{(k)}(p)) \). 
    We estimate the quantitites on the right hand side of 
    \eqref{triangle} in the 
    following way. First of all, since \(  y \in \mathcal 
    E^{(k+1)}(p) \), the distance of the iterates \( \varphi^{j}(y) \) 
    from \( p \) can be calculated by 
  \[
   d(p, \varphi^{j}(y) \leq 
   \int_{p}^{y} 
    \|e^{(k+1)}_{j}(s)\| ds \leq \overline{\| e^{(k+1)}_{j}\|}\eta
    \leq K \bar E_{j} \eta
   \]
      where the integral is taken along \( \mathcal 
    E^{(k+1)}(p) \) and the last inequality follows from Lemma 
    \ref{anglecontraction}. 
   For the second term, the Mean Value Theorem implies
  \[
  d(\varphi^{j}(y),  \varphi^{j}(x)) \leq  \bar F_{j} d(y, x) 
    \leq  \bar F_{j} \omega_{k+1}.
    \]
    Substituting the last two inequalities into \eqref{triangle}, we 
    get 
    \begin{equation}\label{triangle 2} 
    d(p, \varphi^{j}(x)) \leq K \bar E_{j} \eta + 
    \bar F_{j} \omega_{k+1}.
    \end{equation}
    We need to show that this is \( \leq \varepsilon \) for all \( 
    j\leq k+1 \). The first term on the right hand side can be taken \( 
    \leq \varepsilon/2 \) by choosing \( \eta \) sufficiently small. 
    So it is enough to show that 
    \[ 
    \bar F_{j} \omega_{k+1} = K \eta e^{\eta K} \bar H_{k+1} \bar 
    F_{j} \leq \varepsilon /2
    \] for all \( 1\leq j\leq k+1 \). Since we can choose \( \eta \) 
    small, it is then enough to show that there exists the product 
    \( \bar H_{k+1}\bar F_{j} \) is uniformly bounded above for all 
    \( k\geq 1 \) and all \( k+1\geq j \geq 1 \). 
We recall that we can choose the size \( \varepsilon  \) of the 
neighbourhood of \( p \) in which the entire construction takes place 
arbitrarily small, in order to ensure that the constant \( \delta \) 
is arbitrarily small, where \( \delta \) is such that 
    \[ 
K ( \lambda_{u}+ \delta)^{j} \geq F_{j} \geq 
( \lambda_{u}- \delta)^{j} \geq 
(\lambda_{s}+ \delta)^{j} 
\geq  E_{j} \geq 
K^{-1} (\lambda_{s}- \delta)^{j}
\]
This gives 
\[ 
\bar F_{j}\bar H_{j} \leq K (\lambda_{u}+\delta)^{j} 
\frac{(\lambda_{s}+\delta)^{k}}{(\lambda_{s}-\delta)^{k}} 
\leq K 
\left[\frac{(\lambda_{u}+\delta)(\lambda_{s}+\delta)}{\lambda_{u}-\delta}
\right]^{k} \leq K
\]
if \( \delta \) is small enough. 
  \end{proof}

\section{Smoothness}
\label{smoothness}
\begin{lemma}\label{l:smoothness}
    \( \mathcal E^{(\infty)}(z) \) is \( C^{1} \) with 
     Lipschitz continuous derivative. 
    \end{lemma}

\begin{proof}
To see that \( \mathcal E^{(\infty)}(z) \) is \( C^{1} \), we just
substitute  \eqref{uniformpointwiseconvergence} into
\eqref{uniformconvergence} to get
\begin{equation*}
    \|e^{(k)}(z^{(k)}_{t})
    -e^{(k+1)}(z^{(k+1)}_{t}) \| \leq 
     t K \bar H_{k} e^{t K} \leq K\eta e^{\eta K} \bar H_{k}
\end{equation*}
for every \( -\eta \leq t \leq \eta \). 
Since \( \bar H_{k} \to 0 \) exponentially fast, 
this implies that  the sequence
\( e^{(k)}(z_{t}) \) is uniformly Cauchy in \( t \). 
Thus by a standard result about the uniform convergence of derivatives, 
they converge to the tangent directions of the limiting curve 
\( \mathcal E^{(\infty)}(z) \) and this curve is \( C^{1} \).
To prove that the tangent direction are Lipschitz continuous
functions, notice that by the triangle inequality we have
\[ \begin{split}
|e^{\infty}(x)-&e^{\infty}( x' )| \leq \\
& |e^{\infty}(x) - e^{(k)}(x)|  
+
|e^{(k)}(x) - e^{(k)}(x')| + |e^{(k)}(x) - e^{\infty}(x')|
\end{split}
\]
for any two given points \( x, x' \in\mathcal E^{(\infty)}(z)\) and
any \( k\geq k_{0} \).  By Lemma \ref{angleconvergence} we have 
\( |e^{\infty}(x) - e^{(k)}(x)| \leq K \sum_{j\geq k} H_{j}\), 
\( |e^{\infty}(x') - e^{(k)}(x')| \leq K \sum_{j\geq k} H_{j}\) and, by Lemma 
\ref{derivativeconvergence} and the Mean Value Theorem, 
\( |e^{(k)}(x) - e^{(k)}(x')| \leq \|De^{(k)}\| \cdot |x-x'| \leq
 K  |x-x'|\).  This gives 
\[ 
|e^{\infty}(x)-e^{\infty}( x' )| \leq K |x-x'| + K
\sum_{j\geq k} H_{j}.
\]
Since this inequality holds for any \( k \) and \( \sum_{j\geq k}
H_{k} \to 0 \) as \( k\to \infty \) 
it follows that \( |e^{\infty}(x)-e^{\infty}( x' )| 
\leq K |x-x'| \). 
\end{proof}

\section{Contraction}
\label{s:contraction}

Finally we want to show that our limiting curve 
\( \mathcal E^{\infty}(z) \) behaves like a stable manifold in the 
sense that it contracts as \( k\to \infty \). 
Let\( z_{t} \) denote a parametrization by arclength of the leaf 
\( \mathcal E^{\infty}(z) \) \( z_{0}=z \).  

\begin{lemma}\label{contraction}
 There exists a $\tilde\delta>0$, which can be taken arbitrarily small
 with $\delta$, such that
for any \( \eta \geq  t_{2}>t_{1} \geq -\eta \) and  \(
k\geq 1 \) we have
\[ 
\textstyle{|\varphi^{k}(z_{t_{1}})- \varphi^{k}(z_{t_{2}})| \leq K 
(\lambda_{s}+\tilde\delta)^{k} |z_{t_{1}}-z_{t_{2}}|.}
\]
\end{lemma}

\begin{proof}
 Writing \( e^{(\infty)}_{k}=e^{(k)}+ (e^{(\infty)}- e^{(k)}) \) we
 have by the linearity of the derivative
\[ 
\begin{split}
\|\varphi^{k}(z_{t_{1}}) - \varphi^{k}(z_{t_{2}})| 
& =
\int_{t_{1}}^{t_{2}}\|e^{(\infty)}_{k}\| dt \\ 
& = \int_{t_{1}}^{t_{2}}
\|D\varphi^{k} (e^{(k)}) + D\varphi^{k}(e^{(\infty)}- e^{(k)})\| dt.
\end{split}
\]
 Since  \(
\|D\varphi^{k} (e^{(k)})\| = E_{k} \) and  $\|e^{(\infty)}- e^{(n)}\|
\leq K\sum_{j\geq k} H_{j} $ by Lemma \ref{angleconvergence}, we get 
\[ 
\|\varphi^{k}(z_{t_{1}}) - \varphi^{k}(z_{t_{2}})| 
\leq (\bar E_{k} + K \bar F_{k} \sum_{j\geq k} \bar H_{j})  |t_{2}-t_{1}|
\]
Now \( \bar H_{j} \) is an exponentially decreasing sequence and 
therefore \(  \sum_{j\geq k} \bar H_{j} \leq K \bar H_{k}\) and thus, 
reasoning exactly as in the proof of Lemma \ref{smoothness} we have 
\[ 
\bar F_{k} \sum_{j\geq k} \bar H_{j} \leq \bar F_{k} \bar H_{k} \leq 
K\left[\frac{(\lambda_{u}+\delta)(\lambda_{s}+\delta)}{\lambda_{u}-\delta}
\right]^{k} 
\leq K (\lambda_{s}+\tilde\delta)^{k}
\] 
for some \( \tilde\delta > \delta > 0 \) which can be chosen arbitrarily small 
by choosing \( \varepsilon \) small. 
Moreover we also have \( \bar E_{k} \leq (\lambda_{s}+\delta)^{k} 
\leq (\lambda_{s}+\tilde\delta)^{k} \). Finally, using the fact that 
\( |z_{t_2}-z_{t_{1}}| \approx |t_{2}-t_{1}| \) if \( \eta  \) is 
small we conclude the proof. 
\end{proof}

\section{Uniqueness}
\label{uniqueness}
Finally we need to show that the local stable manifold  we have constructed is 
unique in the sense that there is some neighbourhood \( 
B_{\tilde\eta}(p) \) of \( p \) of size \( 
\tilde\eta \) (perhaps smaller than \( \eta \) ) such that 
\( W^{s}_{\eta}(p) \cap B_{\tilde\eta}(p)  \)
is precisely the set of points which stay in this neighbourhood for 
all time., i.e. every other point must leave this neighbourhood at some 
future time. 

\begin{lemma}
    The stable manifold through $p$ is unique in the sense that 
    for there exists some \( \eta > \tilde\eta > 0 \) for which
 \[
 W^{s}_{\eta}(p) \cap B_{\tilde\eta}(p) 
 =\bigcap_{k\geq 1}\mathcal{N}^{k}_{\tilde \eta}(p)
 \]
\end{lemma}
\begin{proof}
 Suppose, by contradiction that there is some point 
 \(  x\in B_{\tilde\eta}(p)  \) which belongs to 
 \( \bigcap_{k\geq 1}\mathcal{N}^{k}_{\tilde \eta}(p) \) but not to 
 \(   W^{s}_{\eta}(p) \). We show that this point must eventually 
 leave \( B_{\tilde\eta}(p) \). 
 
 We need to use here a slightly more 
 refined version of Lemma \ref{regular growth} concerning the 
 hyperbolic stucture in a neighbourhood of the point \( p \). 
 This is also standard in hyperbolic dynamics and we refer the reader 
 to \cites{KatHas94, Rob95} for full details. 
 The property we need is the existence of an invariant expanding 
 conefield in \(  B_{\tilde\eta}(p)   \). This implies that any 
 ``admissible'' curve, i.e. any curve whose tangent directions lie in 
 this expanding conefield will remain admissible for as long as it 
 remains in \(  B_{\tilde\eta}(p)   \) and will grow in length at an 
 exponential rate bounded below by \( \lambda_{u}-\delta \). 
 Moreover the expanding cones do not contain the direction determined 
 by the stable eigenspace \( E^{s}_{p} \) at the fixed point \( p \). 
 By construction the curve \( W^{s}_{\eta} \) is tangent to \( E^{s}_{p} 
 \) at \( p \) and therefore, for \( \eta \) small enough, the tangent 
 vectors to \(  W^{s}_{\eta} \) are not contained in any unstable cone 
 at any point of \(  W^{s}_{\eta} \). Thus \(  W^{s}_{\eta} \) is in 
 some sense \emph{transversal} to the unstable conefield.

 This transversality implies that for any \( x\in  B_{\tilde\eta}(p) 
 \) there exists an admissible curve \( \gamma \) joining \( x \) to some point \( 
 y\in W^{s}_{\eta}  \). The iterates of the curve \( \gamma \) 
 continue to be admissible and the length of \( \gamma \) grows 
 exponentially fast. Therefore, at least one of the endpoints of \( 
 \gamma \) must eventually leave \( B_{\tilde\eta}(p)  \) 
 (or even \( B_{\varepsilon}(p) \)).  By construction, the point \( 
 y\in W^{s}(p) \) never leaves this neighbourhood and therefore the 
 point \( x \) must at some point leave. 
 
\end{proof}

\section{Proof of Lemma \ref{derivativeconvergence}} 

    Since \(D\varphi^{k} \) is a linear map, we have
\begin{equation}\label{hyp0}
     \tan \phi^{(k)} = H_{k+1} \tan \phi^{(k)}_{k+1}. 
\end{equation}
Differentiating on
both sides and taking norms we have
     \begin{equation}\label{hyp3}
	 \begin{split}
	 \|D\phi^{(k)}\| &\leq  \| H_{k+1} \cdot D (\tan
     \phi_{k+1}^{(k)})\| + \| D H_{k} \cdot \tan \phi^{(k)}_{k+1}\| \\ 
     &\leq 
     \|H_{k+1} (1+\tan^{2}\phi^{(k)}_{k+1}) D
     \phi_{k+1}^{(k)}\| + \|D H_{k} \cdot \tan \phi^{(k)}_{k+1} \|
     \end{split}
     \end{equation}
We shall show in three separate sublemmas, 
that \( \phi^{(k)}_{k+1} \leq K\) , \(
D\phi_{k+1}^{(k)} \leq K (F_{k}+E_{k}) \), and \( D H_{k} \leq K
E_{k}\).  Using the fact that \( H_{k}/H_{k+1} \) is also
bounded, and substituting these estimates into \eqref{hyp3} yields 
the estimate in Lemma \ref{derivativeconvergence}. 

\begin{sublemma} \( |\phi^{(k)}_{k+1}| \leq K\). 
    \end{sublemma}
    \begin{proof} 
	Follows immediately by substituting the results of Lemma 
	\ref{angleconvergence} into \eqref{hyp0}
\end{proof}
    \begin{sublemma} \( \|D\phi^{(k)}_{k+1}\| \leq K F_{k}\)
  \end{sublemma}
\begin{proof}
Writing \( \phi^{(k)}_{k+1}=\theta^{(k+1)}_{k+1}-\theta^{(k)}_{k+1} 
    \) we have  
\( \|D\phi^{(k)}_{k+1}\|=\|D\theta^{(k+1)}_{k+1}
    -D\theta^{(k)}_{k+1}\| \leq \|D\theta^{(k+1)}_{k+1}\|
    + \|D\theta^{(k)}_{k+1}\|
 \)
    Our strategy therefore is to obtain estimates for the terms on 
    the right hand side.
First of all we write 
\[
D\varphi^{n}(z)=
\begin{pmatrix}
A_n & B_n \\
C_n & D_n
\end{pmatrix}
\]
where $A_n,B_n,C_n$ and $D_n$ are the matrix entries for the derivative 
$D\varphi^{n}$ evaluated at $z$.
Since $\{e^{(n)}(z),f^{(n)}(z)\}$ correspond
to (resp.) maximal contracting and expanding vectors under $D\varphi^n(z)$
we may obtain them precisely by solving the differential equation
\( \frac{d}{d\t}\| D\varphi^{n}_{z} ( \cos\t, \sin\t) \|=0 \). By
solving this differential equation (and by solving a
similar one for the inverse map $D\varphi^{-n}$) we get  
\[ \tan 2\t^{(k)} =
\frac{2(A_k B_k+C_k D_k)}{A^{2}_{k}+C^{2}_{k}-B^{2}_{k}-D^{2}_{k}}
:=\frac{2\A_k}{\B_k} \]
and 
\[\tan 2\t^{(k)}_{k} = \frac{2(B_k
D_k+A_k C_k)}{D^{2}_{k}+C^{2}_{k}-A^{2}_{k}-B^{2}_{k}}
:=-\frac{2{\mathcal{C}}_k}{{\mathcal{D}}_k}
   \]   
Notice the use of \( \mathcal A_{k}, \mathcal B_{k}, \mathcal C_{k}, \mathcal D_{k} 
\)  as a shorthand notation for the expression in the quotients. 
Now $e^{(k)}_{k},\,f^{(k)}_{k}$ are respectively maximally expanding and
contracting for $D\p^{-k}$, and so we have the identity
\[
D\p^{-k}(\p^{k}(\xi_0))\cdot\textrm{det}\,D\p^{k}(\xi_0)=
\begin{pmatrix}
D_k & -B_k \\
-C_k & A_k
\end{pmatrix}
\]
Then, using the quotient rule for differentiation immediately gives
\begin{equation}\label{anglederivatives}
    \| D\t^{(k)}\| =
\biggl\|\frac{\A_{k}'\B_k-\A_k\B_{k}'}{4\A^{2}_{k}+\B^{2}_{k}}\biggr\|
\text{ and } 
\| D\t^{(k)}_{k}\| =
\biggl\|\frac{\mathcal{D}_{k}'\mathcal{C}_k-\mathcal{D}_k\mathcal{C}_{k}'}
{4\mathcal{C}^{2}_{k}+\mathcal{D}^{2}_{k}}\biggr\|.
\end{equation}

\begin{claim}\label{anglederivatives1}
      \[
      |\mathcal A_{k}|, |\mathcal B_{k}|, |\mathcal C_{k}|,
    |\mathcal D_{k}|  \leq 4\|D\varphi^{k}\|^{2}
\]
and 
\[
    \|\mathcal A'_{k}\|, \|\mathcal B'_{k}\|, |\mathcal C'_{k}\|,
    \|\mathcal D'_{k}\|  \leq 16\|D\varphi^{k}\| \|D^{2}\varphi^{k}\|
    \]
 \end{claim}
\begin{proof}
    For the first set of estimates observe that 
    each of the partial derivatives \( A_{k}, B_{k}, C_{k}, D_{k} \) of 
    \( D\varphi^{k} \) is \( \leq \|D\varphi^{k}\| \). Then 
    \( |\mathcal A_{k}| =  |A_k B_k+C_k D_k| \leq 
    2\|D\varphi^{k}\|^{2}\). The same reasoning gives the estimates 
    in the other cases. 
    To estimate the derivatives, write 
    \(  \|\mathcal A'_{k}\| = |A'_{k}B_{k}+A_{k}B'_{k}+ 
    C'_{k}D_{k}+C_{k}D'_{k}| \).  Now \( |A'_{k}| \leq 2 
    \|D^{2}\varphi^{k}\| \) and similarly for the other terms. 
    \end{proof}

\begin{claim}\label{anglederivatives2}
    \[ 
    4\mathcal{C}^{2}_{k}+\mathcal{D}^{2}_{k}=4\A^{2}_{k}+\B^{2}_{k}=
    (E^{2}_{k}-F^{2}_{k})^2.
    \]
  \end{claim}
  \begin{proof}
      Notice first of all that \( E_{k}^{2}, F_{k}^{2} \) are 
      eigenvalues of  
      \[
      \begin{split}
      (D\p^{k})^{T}D\p^{k} 
      &= 
      \begin{pmatrix}
	  A_{k} & B_{k} \\ C_{k} & D_{k} 
\end{pmatrix}
 \begin{pmatrix}
	  A_{k} & C_{k} \\ B_{k} & D_{k} 
\end{pmatrix}
\\
&= 
 \begin{pmatrix}
	  A^{2}_{k}+B_{k}^{2} & A_{k}C_{k}+ D_{k}B_{k}
	  \\ A_{k}C_{k}+ D_{k}B_{k} & C^{2}_{k} + D^{2}_{k} 
\end{pmatrix}
\end{split}
\]
In particular \( E_{k}^{2}, F_{k}^{2} \) are the two roots of the 
characteristic equation
\( 
\lambda^{2}-\lambda (A_{k}^{2}+B_{k}^{2}+C_{k}^{2}+D_{k}^{2}) + 
 (A^{2}_{k}+B_{k}^{2}) (C^{2}_{k} + D^{2}_{k}) -
(A_{k}C_{k} + B_{k}D_{k})^{2} = 0
\)
and therefore, by the general formula for quadratic equations, we have 
\( 
F_{k}^{2}+E_{k}^{2} = A_{k}^{2}+B_{k}^{2}+C_{k}^{2}+D_{k}^{2}
\) and
\(
E_{k}^{2} F_{k}^{2} =  (A^{2}_{k}+B_{k}^{2}) (C^{2}_{k} + 
D^{2}_{k}) - (A_{k}C_{k} + B_{k}D_{k})^{2}.
\)
From this one can easily check that
\(4\mathcal{C}^{2}_{k}+\mathcal{D}^{2}_{k}=4\A^{2}_{k}+\B^{2}_{k}=
    (E^{2}_{k}-F^{2}_{k})^2=(E^{2}_{k}+F^{2}_{k})^2-4E^{2}_{k}F^{2}_{k}.\)
\end{proof}

Substituting the estimates of Claims
\ref{anglederivatives1}-\ref{anglederivatives2} into
\eqref{anglederivatives} and using the fact that \( k\geq k_{0} \) 
and Lemma \ref{distortion}, this gives
\begin{align*}
    \|D\theta^{(k)}\|, \|D\theta^{(k)}_{k}\| &\leq 128
\frac{\|D\varphi^{k}\|^{3} \|D^{2}\varphi^{k}\|}
{(E_{k}^{2}-F_{k}^{2})^{2}} 
\leq K \frac{\|D\varphi^{k}\|^{3}
\|D^{2}\varphi^{k}\|} {F_{k}^{4}}
\\ &\leq K \frac{
\|D^{2}\varphi^{k}\|}{F_{k}} \leq K F_{k}.
\end{align*}
To estimate $D\t^{(k)}_{k+1}$ we write
$e^{(k)}_{k+1}=\tilde{E}_{k+1}
(\cos\t^{(k)}_{k+1},\sin\t^{(k)}_{k+1})$, so that
\begin{equation*}
\tan\t^{(k)}_{k+1}= \frac{C_1(z_k)\cos\t^{(k)}_{k}+D_1(z_k)\sin\t^{(k)}_{k}}
{A_1(z_k)\cos\t^{(k)}_{k}+B_1(z_k)\sin\t^{(k)}_{k}}=\frac{\M_k}{\N_k}.
\end{equation*}
Notice that we have 
\[ 
\|\mathcal M_{k}\| \leq 2 \|D\varphi_{z_{k}}\| \text{ and } 
\|D\mathcal M_{k}\| \leq 2( \|D^{2}\varphi_{z_{k}}\| \cdot 
\|D\varphi^{k}_{z}\| + \|D\varphi_{z_{k}}\| \cdot 
\|D\theta^{(k)}_{k}\|)
\]
and similarly for \( \|\mathcal N_{k}\| \) and \( \|D\mathcal N_{k}\| \).
Therefore, writing
\begin{equation*}
\|D\t^{(k)}_{k+1}\|= \biggl\|\frac{\N_k\M_{k}'-\M_k\N_{k}'}
{\M^{2}_{k}+\N^{2}_{k}}\biggr\|
\end{equation*}
with
\begin{equation*}
\M^{2}_{k}+\N^{2}_{k}=\frac{\|e^{(k)}_{k+1}\|^2}{\|e^{(k)}_{k}\|^2}\geq
\frac{1}{\|D\p^{-1}(z_k)\|^2},
\end{equation*}
and substituting the estimate obtained,
we get
the desired result. 
\end{proof}
\begin{sublemma}
    \( 
    \|DE_{k}\|, \|DF_{k}\| \leq K F^{2}_k
\quad\text{ and } \quad \|D H_{k}\| \leq K E_k
    \)
 \end{sublemma}
    \begin{proof}
We first estimate
$D_z E_{k}=D\|e^{(k)}_{k}\|$. The corresponding estimate for $D_z F_k$ is
identical. First of all one has,
\(
D_{z} e^{(k)}_{k} =\,D^2\varphi^k(z)e^{(k)}+D\varphi^k\cdot De^{(k)}
\)
and hence
\(
\|D_{z} e^{(k)}_{k}\| \leq\,\|D^2\varphi^k(z)\|+\|D\varphi^k(z)\|\cdot
\|D_{z} e^{(k)}\|.
\)
Since
$D\|e^{(k)}_{k}\|=(e^{(k)}_{k}\cdot De^{(k)}_{k})\|e^{(k)}_{k}\|^{-1}$
one has that 
$\|D_z E_k\|\leq\,\|D e^{(k)}_{k}\|.$ By the previous lemma and the
distortion conditions we get 
\(  \|DE_{k}\|\leq K  F^{2}_{k}\).
 Using the fact that \( \det D\varphi^{k} =
E_{k}F_{k} \) and the quotient rule for differentiation, we get
    \[ 
    DH_{k}=D\left(\frac{E_{k}}{F_{k}}\right) = D\left(\frac{\det
    D\varphi^{k}}{F_{k}^{2}}\right) = \frac{D(\det
    D\varphi^{k})}{F_{k}^{2}} - \frac{2 E_{k} DF_{k}}{F_{k}^{2}}.
    \]
By the estimates above we then get 
\(  \|DH_{k}\| \leq K E_{k}\).
    \end{proof}
This completes the proof of the first inequality in the statement of 
Lemma \ref{derivativeconvergence}. The second one follows immediately 
from the first and the inequality in \eqref{norms}.

\end{document}